\documentclass{tran-l}
\usepackage{amssymb,amscd,amsthm,amsmath,amsfonts}

\newcommand{\B}{\ensuremath{\mathcal B}}
\newcommand{\C}{\ensuremath{\mathbb C}}
\renewcommand{\H}{\ensuremath{\mathcal H}}
\newcommand{\N}{\ensuremath{\mathbb N}}
\newcommand{\R}{\ensuremath{\mathbb R}}

\newcommand{\cs}{{C^{\star}}} 

\newcommand\<{\langle}
\renewcommand\>{\rangle}

\newcommand{\ie}{{\em i.e., }}

\newtheorem{thm}{Theorem}[section] 
\newtheorem{dfn}[thm]{Definition}

\newtheorem{cor}[thm]{Corollary}

\newtheorem{lem}[thm]{Lemma}

\newcommand{\Pf}{{\em Proof}. }

\newtheorem{prop}[thm]{Proposition}

\newtheorem{rmk}[thm]{Remark}
\newtheorem{thm*}{Theorem}
\newcommand{\EPf}{\hfill $\Box$ \medskip}
  {\addvspace{\bigskipamount}\noindent{\bf Example\ }}%

\begin{document}

\title[nonexistence of nontrivial $n$-homomorphisms] {On the nonexistence of
nontrivial involutive $n$-homomorphisms of $\cs$-algebras}
\author{Efton Park and Jody Trout}
\address{Box 298900, Texas Christian University, Fort Worth, TX 76129}
\email{e.park@tcu.edu}
\address{6188 Kemeny Hall, Dartmouth College, Hanover, NH 03755}
\email{jody.trout@dartmouth.edu}
\thanks{MSC $2000$ Classification: Primary 46L05; Secondary 47B99, 47L30} 

\begin{abstract}
An $n$-homomorphism between algebras is a linear map  $\phi : A \to B$ such that
$\phi(a_1 \cdots a_n) = \phi(a_1)\cdots \phi(a_n)$ for all elements $a_1, \dots,
a_n \in A.$ Every homomorphism is an $n$-homomorphism, for all $n \geq 2$, but
the converse is false, in general. Hejazian {\em et al.} \cite{HMM} ask: Is
every $*$-preserving $n$-homomorphism between $\cs$-algebras continuous? We
answer their question in the affirmative, but the even and odd $n$ arguments are
surprisingly disjoint. We then use these results to prove stronger ones: If $n >2$
is even, then $\phi$ is just an ordinary $*$-homomorphism. If $n \geq 3$ is odd, then $\phi$ is a difference of two orthogonal $*$-homomorphisms. Thus, there are no {\it nontrivial} $*$-linear $n$-homomorphisms between $\cs$-algebras.
\end{abstract}

\maketitle

\section{Introduction}

Let $A$ and $B$ be algebras and $n \geq 2$ an integer. A linear map $\phi : A
\to B$ is an {\it $n$-homomorphism} if for all $a_1, a_2,
\dots, a_n \in A$,
\[
\phi(a_1 a_2 \cdots a_n) = \phi(a_1)\phi(a_2) \cdots \phi(a_n).
\]
A $2$-homomorphism is then just a homomorphism, in the usual sense, between
algebras. Furthermore, every homomorphism is clearly also an $n$-homomorphism
for all $n \geq 2$, but the converse is false, in general. The concept of
$n$-homomorphism was studied for complex algebras by Hejazian, Mirzavaziri, and
Moslehian \cite{HMM}. This concept also makes sense for rings and (semi)groups.
For example, an {\it $AE_n$-ring} is a ring $R$ such that every additive
endomorphism $\phi : R \to R$ is an $n$-homomorphism; Feigelstock \cite{Feig1,
Feig2} classified all unital $AE_n$-rings.

In \cite{HMM}, Hejazian {\em et al.} ask: Is every $*$-preserving
$n$-homomorphism between $\cs$-algebras continuous? We answer in the affirmative by proving that every involutive $n$-homomorphism $\phi : A \to B$ between $\cs$-algebras is in fact norm contractive: $\| \phi \| \leq 1$. Surprisingly, the arguments for the even and odd $n$ cases are disjoint and, thus, are discussed in different sections. When $n = 3$, automatic continuity is reported by Bra\v{c}i\v{c} and Moslehian \cite{BM06}, but note that the proof of their Theorem 2.1 does not extend to the nonunital case since the unitization of a $3$-homomorphism is not a $3$-homomorphism, in general.

Using these automatic continuity results, we prove the following stronger
results: If $n > 2$ is even, every $*$-linear $n$-homomorphism $\phi : A \to B$
between $\cs$-algbras is in fact a $*$-homomorphism. If $n \geq 3$ is odd, every
$*$-linear $n$-homomorphism $\phi : A \to B$  is a difference $\phi(a) =
\psi_1(a) - \psi_2(a)$ of  two orthogonal $*$-homomorphisms $\psi_1  \perp
\psi_2$. Regardless, for all integers $n \geq 3$, every {\it positive} linear
$n$-homomorphism is a $*$-homomorphism.   Note that if $\psi$ is a
$*$-homomorphism, then $- \psi = 0 - \psi$ is a norm contractive $*$-preserving
$3$-homomorphism that is not positive linear.

There is also a dichotomy between the unital and nonunital cases.
When the domain algebra $A$ is unital, there is a simple representation of an
$n$-homomorphism as a certain $n$-potent multiple of a homomorphism
(discussed in the Appendix.)  The nonunital case is more subtle. For
example, if $A$ and $B$ are nonunital (Banach) algebras such that $A^n = B^n =
\{0\}$, then {\it every} linear map $L : A \to B$ (bounded or unbounded) is,
trivially, an $n$-homomorphism (see Examples 2.5 and 4.3 of \cite{HMM}). 

The outline of the paper is as follows: In Section $2$, we prove automatic
continuity for the even case and in Section $3$ for the odd case. In Section
$4$, we prove our nonexistence results. A key fact in many of our proofs is the
Cohen Factorization Theorem \cite{Coh59} of $\cs$-algebras. (See Proposition
2.33 \cite{RW98} for an elementary proof of this important result.)  Finally, in Appendix A, we collect some facts about $n$-potents that we need.

The authors would like to thank Dana Williams and Tom Shemanske for their helpful comments and suggestions.

\section{Automatic Continuity: The Even Case}

In this section, we prove that when $n > 2$ is even, every involutive (\ie
$*$-linear) $n$-homomorphism between $\cs$-algebras is completely positive and
norm contractive, which generalizes the well-known result for
$*$-homomorphisms ($n=2$). Recall that a linear map $\theta : A \to B$ between
$\cs$-algebras is {\it positive} if $a \geq 0$ implies $\theta(a) \geq 0$ or,
equivalently, for every $a \in A$ there is a $b \in B$ such that $\theta(a^*a) =
b^*b$. We say that $\theta$ is {\it completely positive} if, for all $k \geq 1$,
the induced map $\theta_k : M_k(A) \to M_k(B)$, $\theta_k((a_{ij})) =
(\theta(a_{ij}))$, on $k \times k$ matrices is positive.
 
\begin{thm} Let $\H$ be a Hilbert space. If $n \geq 2$ is even, then every involutive
$n$-homomorphism from a C*-algebra $A$ into $\B(\H)$ is completely positive.
\end{thm}

\Pf Let $\phi : A \to \B(\H)$ be an involutive $n$-homomorphism. We may assume
$n = 2k > 2$. Let $\<\cdot , \cdot \>$ denote the inner product on
$\H$. By Stinespring's Theorem \cite{Sti55} (see Prop. II.6.6 \cite{Bla06}),
$\phi$ is completely positive if
and only for any $m > 1$ and elements $a_1, \dots, a_m \in A$ and vectors $v_1,
\dots, v_m \in \H$ we have
\[
\sum_{i, j = 1}^m \< \phi(a_i^* a_j)v_j, v_i \> \geq 0.
\]
We proceed as follows: for each $1 \leq i \leq m$ use the Cohen Factorization
Theorem \cite{Coh59} to factor $a_i = a_{i1} \cdots a_{ik}$ into a product of
$k$ elements. Thus, their adjoints factor as $a_i^* = a_{ik}^* \cdots a_{i1}^*$.
Since $n = 2k$, we  compute
\[
\begin{aligned} \sum_{i, j = 1}^m \< \phi(a_i^* a_j)v_j, v_i \> & 
= \sum_{i, j = 1}^m \< \phi(a_{ik}^* \cdots a_{i1}^* a_{j1} \cdots a_{jk}
)v_j, v_i \> \\
& = \sum_{i, j = 1}^m \< \phi(a_{ik})^* \cdots \phi(a_{i1})^* \phi(a_{j1})
\cdots \phi(a_{jk})v_j, v_i \> \\ 
&  = \< \sum_{j=1}^m
\phi(a_{j1}) \cdots \phi(a_{jk})v_j, \sum_{i=1}^m \phi(a_{i1})\cdots\phi(a_{ik})
v_j\> \\
& =\< x, x \> \geq 0, \end{aligned}
\]
where $x = \sum_{i=1}^m \phi(a_{i1})\cdots\phi(a_{ik}) v_i \in \H$.
The result now follows. \EPf

Even though the previous result is a corollary of the more general
theorem below, we have included it because the proof technique is different.

\begin{lem} Let $\phi : A \to B$ be an $n$-homomorphism. Then, for all $k \geq
1$, the induced maps $\phi_k : M_k(A) \to M_k(B)$ on $k \times k$ matrices are
$n$-homomorphisms. Moreover, if $\phi$ is involutive $(\phi(a^*) = \phi(a)^*)$,
then each $\phi_k$ is also involutive.
\end{lem}

\Pf Given $n$ matrices $a^1 = (a^1_{ij}), \dots, a^n = (a^n_{ij})$ in $M_k(A)$,
we can express their product $a^1 a^2 \cdots a^n = (a_{ij})$, where the
$(i, j)$-th entry $a_{ij}$ is given by the formula
\[
a_{ij} = \sum_{m\sb{1}, \cdots, m\sb{n-1} = 1}^k a^1_{i
m\sb1}a^2_{m\sb{1}m\sb{2}} \cdots a^n_{m\sb{n-1}j}.
\]
Since $\phi_k(a^1 a^2 \cdots a^n) = ( \phi(a_{ij}))$ by definition and
\[
\begin{aligned}
\phi(a_{ij}) & = \sum_{m\sb{1}, \cdots, m\sb{n-1} = 1}^k \phi( a^1_{i
m\sb1}a^2_{m\sb{1}m\sb{2}} \cdots a^n_{m\sb{n-1}j}) \\
& = \sum_{m\sb{1}, \cdots, m\sb{n-1}= 1}^k \phi(a^1_{i m\sb1})
\phi(a^2_{m\sb{1}m\sb{2}}) \cdots \phi(a^n_{m\sb{n-1}j}) \\
& = [\phi_k(a^1)\phi_k(a^2) \cdots \phi_k(a^n)]_{ij},
\end{aligned}
\]
it follows that $\phi_k :M_k(A) \to M_k(B)$ is an $n$-homomorphism. 
Now suppose that $\phi$ is involutive.  We compute for all $a = (a_{ij}) \in
M_k(A)$:
\[
\phi_k(a^*) = \phi_k((a_{ji}^*)) = (\phi(a_{ji}^*)) = (\phi(a_{ji})^*) =
\phi_k(a)^*
\]
and hence each $\phi_k :M_k(A) \to M_k(B)$ is involutive. \EPf

\begin{thm}\label{thm:autocty}
Let $\phi : A \to B$ be an involutive $n$-homomorphism between C*-algebras. If
$n \geq 2$ is even, then $\phi$ is completely positive. Thus, $\phi$ is bounded.
\end{thm}

\Pf We may assume $n = 2k > 2$. Since $\phi$ is linear, we want to show that for
every $a \in A$ we have $\phi(a^*a) \geq 0$.  By the Cohen Factorization
Theorem, for any $a \in A$ we can find $a_1, ... , a_k \in A$ such that the
factorization $a = a_1\cdots a_k$ holds. Thus, the adjoint factors as $a^* =
a_k^* \cdots a_1^*$. Since $n = 2k$ and $\phi$ is $n$-multiplicative and
$*$-preserving,
\[
\begin{aligned}
\phi(a^*a) &  = \phi(a_k^* \cdots a_1^* a_1 \cdots a_k) \\
& = \phi(a_k)^* \cdots \phi(a_1)^* \phi(a_1) \cdots \phi(a_k) \\
& = (\phi(a_1) \cdots \phi(a_k))^* (\phi(a_1) \cdots \phi(a_k)) \\
& = b^*b \geq 0,
\end{aligned}
\]
where $b = \phi(a_1) \cdots \phi(a_k) \in B$. Thus, $\phi$ is a positive linear
map. By the previous lemma, all of the induced maps $\phi_k :  M_k(A) \to
M_k(B)$ on $k \times k$ matrices are involutive $n$-homomorphisms and are
positive. Hence, $\phi$ is completely positive and therefore bounded
\cite{Bla06}.
\EPf

We now wish to show that if $n \geq 2$ is even, then an involutive
$n$-homomorphism is actually norm-contractive. First, we will need
generalizations of the familiar $\cs$-identity appropriate for
$n$-homomorphisms.

\begin{lem}\label{lem_identity} Let $A$ be a $\cs$-algebra. For all $k \geq 1$,
we have that
\[
\Bigg\{ \begin{aligned}
& \|x\|^{2k} = \|(x^*x)^k \| \\
& \|x\|^{2k+1} = \| x(x^*x)^k\| 
\end{aligned}
\]
for all $x \in A$.
\end{lem}

\Pf In the even case, we have easily that
\[
\|x\|^{2k} = (\|x\|^2)^k = \|x^*x\|^k = \|(x^*x)^k\|
\]
by the functional calculus since $x^*x \geq 0$. In the odd case, we compute
again using the $\cs$-identity and functional calculus:
\[
\begin{aligned}
\|x(x^*x)^k\|^2 & = \|(x(x^*x)^k)^* (x(x^*x)^k)\| \\
&  = \|(x^*x)^k x^*x(x^*x)^k\| \\
& = \|(x^*x)^{2k+1}\| = \|(x^*x)\|^{2k+1} \\
& = (\|x\|^2)^{2k+1} = (\|x\|^{2k+1})^2;
\end{aligned}
\]
the result follows by taking square roots. \EPf

\begin{thm} Let $\phi : A \to B$ be an involutive $n$-homomorphism of
$\cs$-algebras. If $\phi$ is bounded, then $\phi$  is norm contractive
$($$\|\phi\| \leq 1$$)$.
\end{thm}

\Pf Suppose  $n = 2k$ is even. Then for all $x \in A$ we have
\[
\phi\big((x^*x)^k\big) = \phi(x^*x \cdots x^*x) = \big(\phi(x^*)\phi(x)\big)^k =
\big(\phi(x)^*\phi(x)\big)^k.
\]
Thus by the previous lemma,
\[
\begin{aligned}
\|\phi(x)\|^{n}& = \|\phi(x)\|^{2k} \\
&= \|(\phi(x)^*\phi(x))^k\| = \|\phi((x^*x)^k)\| \\
& \leq \|\phi\| \|(x^*x)^k\| = \|\phi\| \|x\|^{2k} = \|\phi\| \|x\|^n,
\end{aligned}
\]
which implies that $\|\phi\| \leq 1$ by taking $n$-th roots. 

The proof for the odd case $n = 2k + 1$ is similar. \EPf 

\section{Automatic Continuity: The Odd Case}

The positivity methods above do not work when $n$ is odd, since the negation of
a $*$-homomorphism defines an involutive $3$-homomorphism that is (completely)
bounded, but {\bf not} positive. We need the following slight generalization of Lemma 3.5 of Harris \cite{Har81}.

\begin{lem}\label{lemma_Harris} Let $A$ be a $\cs$-algebra and let $\lambda \neq 0$ and $k \geq 1$. If $a \in A$ then $\lambda
\in \sigma((a^*a)^k)$ if and only if there does not exist an element $c \in A$
with 
\begin{equation}\label{eqn:eltc}
c \, (\lambda-(a^*a)^k) = a.
\end{equation}
\end{lem}

\Pf If $\lambda \not\in \sigma((a^*a)^k)$, then $c = a(\lambda-(a^*a)^k)^{-1}
\in A$ satisfies
\[
c\, (\lambda-(a^*a)^k) = a(\lambda-(a^*a)^k)^{-1}(\lambda-(a^*a)^k) = a.
\]
and so (\ref{eqn:eltc}) holds.

On the other hand, if $\lambda \in \sigma((a^*a)^k)$ then, by the commutative
functional calculus, there is a sequence $\{b_m\}_1^\infty$ in the unitization
$A^+$ with $b_m \not\to 0$ but $d_m =_\text{def} (\lambda - (a^*a)^k)b_m \to 0$.
Since $\lambda \neq 0$ we must have
\[
a^*(aa^*)^{k-1}(a b_m) = (a^*a)^kb_m = \lambda b_m - d_m \not \to 0,
\]
which implies $ab_m \not\to 0$. Hence, there does not exist an element $c \in A$
that can satisfy equation (\ref{eqn:eltc}), since this would imply that 
\[
ab_m = c \, (\lambda-(a^*a)^k)b_m \to 0,
\] which is a contradiction.  This proves the lemma. \EPf

We now prove automatic continuity for involutive $n$-homomorphisms of
$\cs$-algebras for all odd values of $n$. Note that we do not assume that $A$ is
unital, nor do we appeal to the {\it unitization} $\phi^+ : A^+ \to B^+$ of $\phi$,
which is {\bf not} an $n$-homomorphism, in general.

\begin{thm}\label{auto3cty} Let $\phi : A \to B$ be an involutive
$n$-homomorphism between $\cs$-algebras. If $n \geq 3$ is odd, then $\|\phi \|
\leq 1$, \ie $\phi$ is norm contractive.
\end{thm}

\Pf Let $n = 2k+1$ where $k \geq 1$.  Given any $a \in A$ and $\lambda > 0$ such
that that $\lambda \not\in \sigma((a^*a)^k)$, there is, by the previous lemma,
an element $c \in A$ such that
\[
a = c\, (\lambda-(a^*a)^k) = (\lambda c - c(a^*a)^k).
\]
Noting that $c(a^*a)^k$ is a product of $2k+1 = n$ elements in $A$, and $\phi$ is a $*$-linear $n$-homomorphism, we compute:
\[
\begin{aligned}
\phi(a) & = \phi(\lambda c - c(a^*a)^k) = \lambda \phi(c) - \phi(c(a^*a)^k)\\
&  = \lambda \phi(c) - \phi(c)(\phi(a)^*\phi(a))^k = \phi(c)(\lambda - (\phi(a)^*\phi(a))^k)
\end{aligned}
\]
which yields that there is an element $\phi(c) \in B$ with:
\[
\phi(c) (\lambda - (\phi(a)^*\phi(a))^k) = \phi(a).
\]
By the previous lemma, we conclude that $\lambda \not\in
\sigma((\phi(a)^*\phi(a))^k)$. Thus, we have shown the following inclusion of
spectra:
\[
\sigma((\phi(a)^*\phi(a))^k) \subseteq \sigma((a^*a)^k) \cup \{0\}.
\] Therefore, by the spectral radius formula \cite[II.1.6.3]{Bla06} and the generalization of the $\cs$-identity in Lemma \ref{lem_identity}, we must deduce that:
\[
\begin{aligned}
\|\phi(a)\|^{2k} & = \|(\phi(a)^*\phi(a))^k\| \\
& = r((\phi(a)^*\phi(a))^k) \leq r((a^*a)^k) \\
& = \|(a^*a)^k\| = \|a\|^{2k},
\end{aligned}
\]
which implies that $\|\phi(a)\| \leq \|a\|$ for all $a \in A$, as desired.
\EPf

Note that the argument in the previous proof does {\it not} work for $n = 2k$ even, since we would need to employ $(a^*a)^{k-1}a$ which is a product of $2k-1 = n-1$ elements as needed, but not self-adjoint, in general. Thus, we could not appeal to the spectral radius formula for self-adjoint elements and Lemma \ref{lemma_Harris} would not apply. Hence, the even and odd $n$ arguments are essentially disjoint.

\section{Nonexistence of Nontrival Involutive $n$-homomorphisms of
$\cs$-algebras}

Our first main result is the nonexistence of nontrivial
$n$-homomorphisms on unital $\cs$-algebras for all $n \geq 3$. We do the unital
case first since it is much simpler to prove and helps to frame the argument for
the nonunital case.

\begin{thm}\label{thm:nonexist} Let $\phi : A \to B$ be an involutive
$n$-homomorphism between the
$\cs$-algebras $A$ and $B$, where $A$ is unital. If $n \geq 2$ is even, then
$\phi$ is a $*$-homomorphism. If $n \geq 3$ is odd, then $\phi$ is the
difference $\phi(a) = \psi_1(a) - \psi_2(a)$ of two orthogonal $*$-homomorphisms
$\psi_1 \perp \psi_2 : A \to B$.
\end{thm}

\Pf In either case, by Proposition \ref{prop:unitalnhom}, the element $e =
\phi(1) \in B$ is an $n$-potent ($e^n = e$) and is self-adjoint, because
\[
e = \phi(1) = \phi(1^*) = \phi(1)^* = e^*.
\]
Also, there is an associated algebra homomorphism $\psi : A \to B$ defined for
all $a \in A$ by the formula
\[
\psi(a) = e^{n-2} \phi(a) = \phi(a)e^{n-2}
\]
such that $\phi(a) = e\psi(a) = \psi(a) e$. In either case, $\psi$ is $*$-linear since $\phi$ is
$*$-linear and $e$
is self-adjoint and commutes with the range of $\phi$:
\[
\psi(a^*) = e^{n-2} \phi(a^*) = e^{n-2} \phi(a)^* = \big(e^{n-2} \phi(a)
\big)^* = \psi(a)^*.
\]
Now, if $n = 2k$ is even, $e = e^n = (e^k)^* e^k \geq 0$ and so $e = p$ is a
projection. Thus, $\phi(a) = p \psi(a) = \psi(a) p = p \psi(a) p$ is a
$*$-homomorphism. If $n \geq 3$ is odd, then by Lemma \ref{lem_n-potents}, $e$
is the
difference of two orthogonal projections $e = p_1 - p_2$ which must commute with
both $\psi$ and $\phi$ by the functional calculus. Define $\psi_1, \psi_2 : A
\to B$ by $\psi_i(a) = p_i \psi(a) p_i$ for all $a \in A$ and $i = 1, 2$. Then
$\psi_2 \perp \psi_2$ are orthogonal $*$-homomorphisms, and 
\[
\psi_1(a) - \psi_2(a) = p_1 \psi(a) - p_2 \psi(a) = e \psi(a) = \phi(a)
\]
for all $a \in A$, from which the desired result follows.\EPf

\begin{cor} Let $\phi : A \to B$ be a linear map between $\cs$-algebras. If $A$
is unital, the following are equivalent for all integers $n \geq 2$:
\begin{itemize}
\item[a.)] $\phi$ is a $*$-homomorphism.
\item[b.)] $\phi$ is a positive $n$-homomorphism.
\item[c.)] $\phi$ is an involutive $n$-homomorphism and $\phi(1) \geq 0$.
\end{itemize}
\end{cor}
\Pf Clearly (a) $\implies$ (b) $\implies$ (c). If $n \geq 2$ is even, then (c)
$\implies$ (a) by the previous result. If $n \geq 3$ is odd, then by the
previous result, we only need to show that $\phi$ is positive. Let $n = 2k+1$.
Given any $a \in A$, by the Cohen Factorization Theorem, we can write $a = a_1
\cdots a_k$. Since $\phi(1) \geq 0$, by hypothesis, and $n = 2k+1$, we compute:
\[
\begin{aligned}
\phi(a^* a) & = \phi(a^* 1 a) = \phi(a_k^* \cdots a_1^* 1 a_1 \cdots a_k) \\
& = \phi(a_k)^* \cdots \phi(a_1)^* \phi(1) \phi(a_1) \cdots \phi(a_k)  \\
& = \big( \phi(a_1) \cdots \phi(a_k)\big)^* \phi(1) \big(\phi(a_1) \cdots
\phi(a_k)\big) \\
& = b^* \phi(1) b \geq 0,
\end{aligned}
\]
where $b = \phi(a_1) \cdots \phi(a_k) \in B$. Thus, $\phi$ is positive linear
and therefore a $*$-homomorphism. \EPf

Next, we extend our nonexistence results to the nonunital case, by appealing to
approximate unit arguments (which require continuity!) and the following
important factorization
property of $*$-preserving $n$-homomorphisms.

\begin{lem}[Coherent Factorization Lemma] Let $\phi : A \to B$ be an involutive
$n$-homomorphism of $\cs$-algebras. For any $1 \leq k \leq n$ and any $a \in A$,
if $a = a_1 \cdots a_k = b_1 \cdots b_k$ in $A$, then
\[
\phi(a_1) \cdots \phi(a_k) = \phi(b_1) \cdots \phi(b_k) \in B.
\]
\end{lem}

Note that, in general, $\phi(a) \neq \phi(a_1) \cdots \phi(a_k)$ when $1 < k <
n$.

\medskip

\Pf Clearly, we may assume $1 < k < n$. Since $\phi$ is $*$-linear, the range
$\phi(A) \subset B$ is a self-adjoint linear subspace of $B$ (but not
necessarily a subalgebra, in general).  Given any $d = \phi(c) \in \phi(A)$,
using the Cohen Factorization Theorem, write $d = d_1 \cdots d_n = \phi(c_1)
\cdots \phi(c_n)$ where $d_i = \phi(c_i)$ for $1 \leq i \leq n$. Consider the
following computation:
\[
\begin{aligned}
\phi(a_1)\cdots\phi(a_k)d & = \phi(a_1) \cdots \phi(a_k) \phi(c_1) \cdots
\phi(c_n) \\
& = \phi(a_1 \cdots a_k c_1 \cdots c_{n-k}) \phi(c_{n-k+1}) \cdots \phi(c_n)\\
& = \phi(b_1 \cdots b_k c_1 \cdots c_{n-k}) \phi(c_{n-k+1}) \cdots \phi(c_n) \\
& = \phi(b_1) \cdots \phi(b_k) \phi(c_1) \cdots \phi(c_n) \\
& = \phi(b_1) \cdots \phi(b_k) d.
\end{aligned}
\]
Let $f = \phi(a_1) \cdots \phi(a_k) - \phi(b_1) \cdots \phi(b_k)$. Then
$f d =0$ for all $d \in \phi(A) \subset B$, and thus $f d =0$
for all $d$ in the $*$-subalgebra $A_\phi$ of $B$
generated by $\phi(A)$. In particular, for the element
\[
d_a = \phi(a_k^*)\cdots \phi(a_1^*) - \phi(b_k^*) \cdots \phi(b_1^*) = 
f^*\in A_\phi.
\]
Hence, $ff^* = fd_a = 0$ and so $\|f\|^2 = \|ff^*\| = 0$ by the $\cs$-identity.
Therefore,
\[
\phi(a_1) \cdots \phi(a_k) - \phi(b_1) \cdots \phi(b_k) = f = 0,
\]
and the result is proven. \EPf

\begin{dfn} An approximate unit for a $($nonunital$)$ $\cs$-algebra $A$ is a net
$\{e_\lambda\}_{\lambda \in \Lambda}$of elements in $A$ indexed by a directed
set $\Lambda$ such that
\begin{itemize}
\item[a.)] $0 \leq e_\lambda$ and $\|e_\lambda\| \leq 1$ for all $\lambda \in
\Lambda$;
\item[b.)] $e_\lambda \leq e_\mu$ if $\lambda \leq \mu$ in $\Lambda$;
\item[c.)] For all $a \in A$,
\[
\lim_{\lambda \to \infty} \|ae_\lambda - a\| = \lim_{\lambda \to \infty}
\|e_\lambda a - a \| = 0.
\]
\end{itemize}
\end{dfn}
\noindent Every $\cs$-algebra has an approximate unit, which is
countable ($\Lambda = \N$) if $A$ is separable (see Section II.4 of Blackadar
\cite{Bla06}.)

\begin{thm}\label{thm:asociated_hom} Suppose $\phi : A \to B$ is an involutive
$n$-homomorphism of
$\cs$-algebras, where $A$ is nonunital. Then, for all $a \in A$, the limit
\[
\psi(a) = \lim_{\lambda \to \infty} \phi(e_\lambda)^{n-2} \phi(a) =
\lim_{\lambda \to \infty} \phi(a)\phi(e_\lambda)^{n-2}
\]
exists, independently of the choice of the approximate unit $\{e_\lambda\}$ of
$A$, and defines a $*$-homomorphism $\psi : A\to B$ such that
\[
\phi(a) = \lim_{\lambda \to \infty} \phi(e_\lambda) \psi(a)
\]
for all $a \in A$. 
\end{thm}

\Pf We may assume $n \geq 3$. Given $a \in A$, use the Cohen Factorization
Theorem to factor $a = a_1 a_2 \cdots a_n$. Define a map $\psi : A \to B$
by
\[
\psi(a) = \phi(a_1a_2) \phi(a_3) \cdots \phi(a_n) = \phi(a_1) \cdots
\phi(a_{n-2})\phi(a_{n-1}a_n),
\]
which is {\it well-defined} by the Coherent Factorization Lemma.
The continuity of $\phi$ implies that
\[
\begin{aligned}
\lim_{\lambda\to\infty}\phi(e_\lambda)^{n-2} \phi(a) & =
\lim_{\lambda\to\infty}\phi(e_\lambda)^{n-2} \phi(a_1) \cdots \phi(a_n) \\
&= \lim_{\lambda\to\infty}\phi(e_\lambda^{n-2} a_1 a_2) \phi(a_3) \cdots
\phi(a_n) \\
&= \phi(a_1 a_2) \phi(a_3) \cdots \phi(a_n)  = \psi(a)\in B.
\end{aligned}
\]
It follows that we can write:
\[
\psi(a) = \lim_{\lambda \to \infty} \phi(e_\lambda)^{n-2} \phi(a) =
\lim_{\lambda \to \infty} \phi(a)\phi(e_\lambda)^{n-2},
\]
and so $\psi : A \to B$ is linear since $\phi$ is linear. Moreover, since $\phi$
is $*$-linear, it follows that $\psi$ is also $*$-linear: 
\[
\begin{aligned}
\psi(a)^* & = \big(\phi(a_1a_2)\phi(a_3) \cdots \phi(a_n)\big)^* \\
               & = \phi(a_n)^* \cdots \phi(a_3)^* \phi(a_1a_2)^* \\
               & = \phi(a_n^*) \cdots \phi(a_3^*) \phi(a_2^* a_1^*) \\
               & = \phi(a_{n1}^* a_{n2}^*) \phi(a_{n-1}^*) \cdots
\phi(a_{12}^*)\\
               & = \psi((a_{n1}^* a_{n2}^*)(a_{n-1}^*) \cdots (a_{12}^*)) \\
               & = \psi(a_n^* \cdots a_1^*) = \psi(a^*).
\end{aligned}
\]
In the computation above, we factored $a_n = a_{n2} a_{n1}$ and set
$a_{12} = a_1a_2$ to obtain the factorization $a^* = a_n^* \cdots a_1^* =
(a_{n1}^* a_{n2}^*) a_{n-1}^* \cdots
a_{12}^*$ into $n$ elements. Given $a, b \in A$ with factorizations $a = a_1
\cdots a_n$ and $b = b_1 \cdots b_n$, the fact that $\phi$ is
an $n$-homomorphism implies:
\[
\begin{aligned}
\psi(a)\psi(b) & = \big(\phi(a_1a_2) \phi(a_3) \cdots
\phi(a_n)\big)\big(\phi(b_1b_2) \phi(b_3) \cdots \phi(b_n)\big) \\
& = \phi((a_1a_2) a_3 \cdots a_n (b_1b_2)) \phi(b_3) \cdots \phi(b_n) \\
& = \phi((ab_1)b_2) \phi(b_3) \cdots \phi(b_n) \\
& = \psi(ab);
\end{aligned}
\]
note that $ab = (ab_1) b_2 b_3 \cdots b_n$ is a factorization of $ab$ into $n$
elements. A second proof of multiplicativity goes as follows:
\[
\begin{aligned}
\psi(ab) & = \lim_{\lambda \to \infty} \phi(e_\lambda)^{n-2} \phi(ab)
  = \lim_{\lambda \to \infty} \phi(e_\lambda)^{n-2} \phi(\lim_{\mu \to \infty} a
e_\mu^{n-2} b) \\
& = \lim_{\lambda \to \infty} \phi(e_\lambda)^{n-2} \lim_{\mu \to \infty}\phi(a
e_\mu^{n-2} b) \\ 
& = \lim_{\lambda \to \infty} \phi(e_\lambda)^{n-2} \lim_{\mu \to \infty}\phi(a)
\phi(e_\mu)^{n-2} \phi(b) \\
& = \lim_{\lambda \to \infty}\phi(e_\lambda)^{n-2} \phi(a) \lim_{\mu \to
\infty}\phi(e_\mu)^{n-2} \phi(b) \\
& = \psi(a) \psi(b).
\end{aligned}
\]
Thus, $\psi$ is a well-defined $*$-homomorphism. Finally, we compute:
\[
\begin{aligned}
\lim_{\lambda \to \infty}\phi(e_\lambda) \psi(a) &= 
\lim_{\lambda \to \infty}\phi(e_\lambda) \phi(a_1a_2) \phi(a_3) \cdots \phi(a_n)
\\
&= \lim_{\lambda \to \infty}\phi(e_\lambda (a_1a_2) a_3 \cdots a_n) 
=\lim_{\lambda \to \infty} \phi(e_\lambda a) \\
&= \phi(a).
\end{aligned}
\]
\EPf

Using similar factorizations, the fact that $\{e_\lambda^n\}$ is also an
approximate unit for $A$, and the fact that the strict completion of the
$\cs$-algebra $\cs(\phi(A))$ generated by the range $\phi(A)$ is the multiplier
algebra $M(\cs(\psi(A)))$, we obtain the nonunital version of Proposition
\ref{prop:unitalnhom}.

\begin{cor} Suppose that $A$ and $B$ are $\cs$-algebras with $A$ nonunital,
and let $\phi : A \to B$ be an involutive $n$-homomorphism
with associated $*$-homomorphism $\psi: A \to B$. Then there is a 
self-adjoint $n$-potent $e = e^* = e^n \in M(C^*(\phi(A)))$
such that $\phi(e_\lambda) \to e$ strictly for any approximate unit
$\{e_\lambda\}$ of $A$, and with the property that
\[
\begin{aligned}
\phi(a) & = e \psi(a) = \psi(a) e \\
\psi(a) & = e^{n-2} \phi(a) 
\end{aligned}
\]
for all $a \in A$.
\end{cor}

\Pf  By the previous proof, we can define $e \in M(\cs(\phi(A)))$ on generators
$\phi(a)$ by
\[e \phi(a) = \lim_{\lambda \to \infty} \phi(e_\lambda)\phi(a) = \phi(a_1 a_2
\cdots a_{n-1})\phi(a_n) \in \cs(\phi(A))
\]
for any $a = a_1 \cdots a_n \in A$. It follows that:
\[
\begin{aligned}
e^n \phi(a) & = \lim_{\lambda \to \infty} \phi(e_\lambda)^n \phi(a) \\
& = \lim_{\lambda \to \infty} \phi(e_\lambda^n) \phi(a_1) \phi(a_2) \cdots
\phi(a_n) \\
& = \lim_{\lambda \to \infty} \phi((e_\lambda^n) a_1a_2 \cdots a_{n-1})\phi(a_n)
\\
& = \phi(a_1 \cdots a_{n-1}) \phi(a_n) = e \phi(a),
\end{aligned}
\]
which implies $e \in M(\cs(\phi(A)))$ is $n$-potent. The fact that $e = e^*$
follows from $\phi(e_\lambda)^* = \phi(e_\lambda^*) = \phi(e_\lambda)$. The
other statements follow from the previous proof. \EPf

The dichotomy between the unital and nonunital cases is now clear. If $A$ is
unital, then $\cs(\phi(A)) \subset B$ is a unital $\cs$-subalgebra of $B$ with
unit
$\psi(1) = \phi(1)^{n-1} \in B$ (which is a projection!) and so 
\[ M(\cs(\psi(A)))
= \cs(\phi(A)) \subset B.\]
However, for $A$ nonunital, we cannot identify the multiplier algebra
$M(\cs(\phi(A)))$ as a subalgebra of $B$, or even $M(B)$, unless $\phi$ is
surjective. In general, we only have inclusions $\psi(A) \subset \cs(\phi(A))
\subset B.$

Now that we know, as in the unital case, every involutive
$n$-homomorphism is an $n$-potent multiple of a $*$-homomorphism, we can prove
the following general version of Theorem \ref{thm:nonexist} and its corollary in
a similar manner using Lemma \ref{lem_n-potents}.

\begin{thm} Let $\phi : A \to B$ be an involutive $n$-homomorphism of
$\cs$-algebras. If $n \geq 2$ is even, then $\phi$ is
a $*$-homomorphism. If $n \geq 3$ is odd,
then $\phi$ is the difference $\phi(a) = \psi_1(a) - \psi_2(a)$ of two
orthogonal
$*$-homomorphisms $\psi_1 \perp \psi_2 : A \to B$.
\end{thm}

\begin{cor} For all $n \geq 2$ and $\cs$-algebras $A$ and $B$, $\phi : A \to B$
is a positive $n$-homomorphism if and only if $\phi$ is a $*$-homomorphism.
\end{cor}

\appendix
\section{On $n$-homomorphisms and $n$-potents}

An element $x  \in A$ is called an {\it $n$-potent} if $x^n = x$. Note that if
$\phi : A \to B$ is an $n$-homomorphism, then $\phi(x) = \phi(x^n) = \phi(x)^n
\in B$ is also an $n$-potent. The following important result is Proposition 2.2
\cite{HMM}, whose proof is included for completeness.

\begin{prop}\label{prop:unitalnhom}
If $A$ is a unital algebra $($or ring$)$ and $\phi : A \to B$ is an
$n$-homomorphism, then there is a homomorphism $\psi : A \to B$ and an
$n$-potent $e = e^n \in B$ such that $\phi(a) = e \psi(a) = \psi(a) e$ for all
$a \in A$. Also, $e$ commutes with the range\footnote{Note that the range
$\phi(A)$ is not a subalgebra of $B$ in general.} of $\phi$, \ie $e \phi(a) =
\phi(a)
e$ for all $a \in A$.\end{prop}

\Pf Note that $e = \phi(1) = \phi(1^n) = \phi(1)^n = e^n \in B$ is an
$n$-potent. Define a linear map $\psi :A \to B$ by $\psi(a) = e^{n-1} \phi(a)$
for all $a \in A$. For all $a, b \in A$,
\[
\begin{aligned} \psi(ab) & = e^{n-2} \phi(ab) = e^{n-2} \phi(a 1^{n-2} b) \\
& = \big(e^{n-2}\phi(a)\big)\big(\phi(1)^{n-2} \phi(b)\big) \\
& = \big(e^{n-2}\phi(a)\big) \big(e^{n-2}\phi(b)\big)= \psi(a) \psi(b),
\end{aligned}
\]
and so $\psi$ is an algebra homomorphism. Furthermore,
\[
e \psi(a) = \phi(1) (\phi(1)^{n-2}\phi(a)) = \phi(1)^{n-1}\phi(a) =
\phi(1^{n-1} a) = \phi(a).
\]
Similarly, $\psi(a)e = \phi(a)$ for all $a \in A$. The final statement is a
consequence of the fact that for all $a \in A$,
\[
e \phi(a) = \phi(1) \phi(a1^{n-1}) = \big(\phi(1) \phi(a)\phi(1)^{n-2}\big)
\phi(1) = \phi(1a1^{n-2})e =  \phi(a)e.
\]
\EPf

The following computation will be more significant when we consider
the nonunital case (see the proof of Theorem \ref{thm:asociated_hom}.)

\begin{cor} Let $\phi$ and $\psi$ be as in Proposition \ref{prop:unitalnhom} and
$n \geq 3$.  Then for all $a \in A$, if $a = a_1a_2 \cdots a_n$ with $a_1,
\dots, a_n \in A$,
\[
\psi(a) = \phi(a_1a_2) \phi(a_3) \cdots \phi(a_n).
\]
\end{cor}

\Pf We compute as follows:
\[
\begin{aligned}
\psi(a) & =_\text{def} e^{n-2}\phi(a) = \phi(1)^{n-2}\phi(a_1 \cdots a_n) \\
&= \phi(1)^{n-2} \phi(a_1) \cdots \phi(a_n) \\
& = \big(\phi(1)^{n-2}\phi(a_1)\phi(a_2)\big) \phi(a_3) \cdots \phi(a_n) \\
&= \phi(1^{n-2}a_1a_2) \phi(a_3) \cdots \phi(a_n) \\
& = \phi(a_1a_2)\phi(a_3) \cdots \phi(a_n). \hskip1.5in \Box
\end{aligned}
\]

\begin{dfn} Let $A$ be a unital algebra. An {\em $n$-partition of unity} is an
ordered $n$-tuple $(e_0, e_1, \dots, e_{n-1})$ of idempotents $(e_k^2 = e_k)$
that sum to the identity $e_0 + e_1 + \cdots + e_{n-1} = 1$
and are pairwise mutually orthogonal, \ie $e_j e_k = \delta_{jk} 1$
for all $0 \leq j, k \leq n-1$, where $\delta_{jk}$ is the Kronecker delta.
\end{dfn}

Note that $e_0 = 1 - (e_1 + \cdots + e_{n-1})$ is completely determined by $e_1, 
e_2, \dots, e_{n-1}$ and is thus redundant in the notation for an $n$-partition
of unity.

\begin{dfn}\label{def:roots} Let $\omega_ 0 = 0$ and $\omega_k = e^{2 \pi i
(k-1) / ( n-1)}$ for $ 1\leq k \leq n-1$. Note that $\omega_1 = 1$ and
$\omega_1, \dots, \omega_{n-1}$ are the $(n-1)$-th roots of unity and
$\Sigma_n = \{\omega_0, \omega_1, \dots, \omega_{n-1}\}$
are the $n$ roots of the polynomial equation $x^n - x = x(x^{n-1} - 1) = 0$.
\end{dfn} 

If $A$ is a complex algebra, we let $\tilde{A}$ denote $A$, if $A$ is unital, or
the unitization $A^+ = A \oplus \C$, if $A$ is nonunital.

\begin{thm}\label{thm:partition} Let $A$ be a complex algebra. If $e \in A$ is
an $n$-potent, there is a unique $n$-partition of unity $(e_0, e_1, \dots ,
e_{n-1})$ in $\tilde{A}$ such that
\[
e = \sum_{k = 1}^{n-1} \omega_k e_k.
\]
If $A$ is nonunital, then $e_1, \dots, e_{n-1} \in A$.
\end{thm}

\Pf Define the $n$ polynomials $p_0, p_1, \dots, p_{n-1}$ by
\[
p_k(x) = \frac{\prod_{j \neq k} (x - \omega_j)}
{\prod_{j \neq k} (\omega_k - \omega_j)}.
\]
In particular, $p_0(x) = 1 - x^{n-1}$.  Each polynomial $p_k$ has degree $n-1$
and satisfies $p_k(\omega_k) = 1$ and $p_k(\omega_j) = 0$ for all $j \neq k$. It
follows that $p_j(x)p_k(x) = 0$ for all $x \in \Sigma_n$.
We also claim that for all $x \in \C$ that 
\begin{equation}\label{eqn:sumone}
\sum_{k = 0}^{n-1} p_k(x) = p_0(x) + \cdots + p_{n-1}(x) = 1
\end{equation}
\begin{equation}\label{eqn:identity}
x = \sum_{k = 0}^{n-1} \omega_k p_k(x).
\end{equation}
Indeed, these identities follow from the fact that these polynomial equations
have degree $n-1$ but are satisfied by the $n$ {\it
distinct} points in $\Sigma_n$.

Now, given any $x^n = x$ in $\C$ it follows that $p_k(x)^2 = p_k(x)$. Hence, for
any $n$-potent $e \in A$, if we define $e_k = p_k(e)$ then 
$(e_0, e_1, \dots, e_{n-1})$ consists of idempotents $e_k^2 = p_k(e)^2 = p_k(e)
= e_k$ and satisfy, by (\ref{eqn:sumone}),
\[
\sum_{k=0}^{n-1} e_k = \sum_{k=0}^{n-1} p_k(e) = 1_{\tilde{A}}.
\]
They are pairwise orthogonal, because 
$e_j e_k = p_j(e)p_k(e) = 0$ for $j \neq k$.  Moreover,
\[
e = \sum_{k = 1}^{n-1} \omega_k p_k(e) = \sum_{k =1}^{n-1} \omega_k e_k
\]
by Equation (\ref{eqn:identity}). For $1 \leq k \leq n-1$, note that $p_k(x) = x
q_k(x)$ for some polynomial $q_k(x)$. Hence, if $A$ is nonunital and $1 \leq k
\leq n-1$, we have $e_k = p_k(e) = e q_k(e) \in A$, since $A$ is an ideal in
$\tilde{A}$. \EPf

The following result is the $n$-homomorphism version of the previous $n$-potent
result. Recall say that two linear maps $\psi_i, \psi_j : A\to B$ are 
{\it orthogonal} ($\psi_i \perp \psi_j$) if $$\psi_i(a) \psi_j(b) =
\psi_j(b)\psi_i(a)
= 0$$ for all $a, b \in A$.\footnote{Note that the zero homomorphism is
orthogonal to every homomorphism.}

\begin{prop} Let $A$ and $B$ be complex algebras. If $A$ is unital then a linear
map $\phi : A \to B$ is an $n$-homomorphism if and only if there exist $n-1$
mutually orthogonal homomorphisms $\psi_1, \dots , \psi_{n-1} : A \to B$ such
that for all $a \in A$,
\[\phi(a) = \sum_{k =1}^{n-1} \omega_k \psi_k(a).
\]
\end{prop}

\Pf $(\Rightarrow)$ Let $\phi : A \to B$ be an $n$-homomorphism. By Proposition
\ref{prop:unitalnhom}, there is an $n$-potent $e \in B$ and a homomorphism $\psi
: A \to B$  such that $\phi(a) = e \psi(a) = \psi(a) e$. Using the previous
result, write $e = \sum_{k=1}^{n-1} \omega_k e_k$, where $(e_0, e_1, \dots,
e_{n-1})$ is the associated $n$-partition of unity in $\tilde{A}$ defined by the
polynomials
$p_k$. Since $e_k = p_k(e)$, we have that $e_k \psi(a) = \psi(a) e_k$ for $1
\leq k \leq n-1$. Define $\psi_k : A \to B$  by 
\[
\psi_k(a) =_\text{def} e_k
\psi(a) = e_k^2 \psi(a) = e_k \psi(a) e_k.
\]
Then $\psi_1, \dots, \psi_{n-1}$
are orthogonal homomorphisms and, for all $a \in A$,
\[
\phi(a) = e \psi(a) = \sum_{k=1}^{n-1} \omega_k e_k \psi(a) = \sum_{k=1}^{n-1}
\omega_k \psi_k(a).
\]

$(\Leftarrow)$ Follows from the fact that $\omega_k^n = \omega_k$ for all $k =
1, \dots, n-1$. \EPf

\begin{rmk} If $A$ is nonunital, the above result does not hold. One reason is
that the unitization $\phi^+ : A^+ \to B^+$ of an $n$-homomorphism is not, in
general, an $n$-homomorphism. Also, if $A^n = B^n = \{0\}$, 
then every linear map $L : A \to B$ is an $n$-homomorphism
$($See Examples 2.5 and 4.3 of Hejazian et al \cite{HMM}$)$. \end{rmk}

Let $\Sigma_n$ be the $n$ roots of the polynomial equation $x = x^n$ from
Definition \ref{def:roots}. If $A$ is a $\cs$-algebra, it follows that a normal
$n$-potent $e = e^n$ must have spectrum $\sigma(e) \subseteq \Sigma_n$. Recall
that a projection is an element $p = p^* = p^2 \in A$. Two projections $p_1$ and
$p_2$ are orthogonal if $p_1p_2 = 0$. A tripotent is a $3$-potent element $e^3 =
e \in A$.

The following characterization of self-adjoint $n$-potents in $\cs$-algebras is important for our nonexistence results on $n$-homomorphisms.

\begin{lem}\label{lem_n-potents} Let $A$ be a $\cs$-algebra.
\begin{itemize}
\item[a.)] If $n \geq 2$ is an even integer, the following are equivalent:
\begin{itemize}
\item[i.)] $e$ is a projection.
\item[ii.)] $e$ is a positive $n$-potent.
\item[iii.)] $e$ is a self-adjoint $n$-potent.
\end{itemize}
\item[b.)] If $n \geq 3$ is an odd integer, the following are equivalent:
\begin{itemize}
\item[i.)] $e$ is a self-adjoint tripotent.
\item[ii.)] $e = p_1 - p_2$ is a difference of two orthogonal projections.
\item[iii.)] $e$ is a self-adjoint $n$-potent.
\end{itemize}
\end{itemize}
\end{lem}

\Pf In both the even and odd cases,  (i) $\implies$ (ii) $\implies$ (iii) (See
Theorem \ref{thm:partition}).  Suppose (iii) holds. If $n = 2k$ is even, 
\[
e =e^* = e^n = e^{2k} = (e^k)^* (e^{k}) \geq 0,
\]
and so the spectrum of $e$ satisfies $\sigma(e) \subset \Sigma_n \cap [0,
\infty] =
\{0, 1\}$. Thus, $e$ is a projection. If $n \geq 3$ is odd, then since $e = e^*$
we must have $\sigma(e) \subset \Sigma_n \cap \R = \{-1, 0, 1\}$. Thus, $\lambda
= \lambda^3$
for all $\lambda \in \sigma(e)$, which implies $e = e^3$ is tripotent. \EPf


\end{document}